\documentclass[]{article}

\usepackage{amssymb,amsmath, graphics}
\usepackage{color}
\newtheorem{thm}{Theorem}
\newtheorem{prop}{Proposition}
\newtheorem{lem}{Lemma}
\newtheorem{rem}{Remark}
\newtheorem{definition}{Definition}

\begin{document}

\title{{\itshape Testing means from sampling populations \\ with undefined labels}}

\author{Florent AUTIN\footnote{Address : C.M.I., 39 rue F. Joliot Curie, 13453 Marseille Cedex 13. Aix-Marseille Universit\'e. FRANCE. Email: autin@cmi.univ-mrs.fr}\hspace{2pt}
(Aix-Marseille Universit\'e) \\
and Christophe POUET\footnote{Address : Ecole Centrale Marseille, 38 rue F. Joliot-Curie, 13451 Marseille Cedex 20. FRANCE. Email : cpouet@centrale-marseille.fr}
\hspace{2pt} (Ecole Centrale Marseille)}

\date{}

\maketitle

\begin{abstract}
We consider the problem of testing means from samples of two populations  for which the labels are not defined with certainty. We  show that this problem is connected to another one that is testing expected values of components of mixture-models from two data samples.  The underlying mixture-model is associated with known varying mixing-weights.  We provide a testing procedure that performs well. Then we point out  the loss of performance of our method due to the mixing-effect by comparing its numerical performances to the Welch's t-test on means which would have been done if true labels were available.
 \end{abstract}

\noindent {\it{Keywords}}: Asymptotic distribution, hypothesis testing, missing labels, mixture-models.\\

\noindent {\it{AMS Subject Classification}}: 62D05, 62F03, 62F05.

\section{Introduction}

\noindent In many cases, researchers can be interested in gathering information about two populations in order to compare them. In that setting, tests of significance are useful statistical tools for detecting a difference between two population parameters. Related application fields are numerous. Some examples are genetics, neuronal data analysis, medicine, biology, physics, chemistry, social sciences, among other fields.  \\

\noindent Consider the Gaussian setting, for which each data of the two populations under study is assumed to follow a normal distribution. Let us recall that this assumption can be tested beforehand using a normality test, such as the well-known Shapiro-Wilk or Kolmogorov-Smirnov test, or it can be assessed graphically using a normal quantile plot. Comparisons between the means of the populations are usually carried out by using  t-statistic and lead to the well known Student's t-test or Welch's t-test (see Welch \cite{Welch}).\\  

 \noindent These t-tests are popular because of their ease of use and their good performances. Moreover they are robust in the sense that
 they still perform well when the components are not really Gaussian, provided that the samples size are large enough. Nevertheless, these
 testing methods require to know the {\it{label}} of each data, that is the population each data is associated with. Unfortunately, researchers do
 sometimes not get this information. Indeed, one can imagine some cases where the labels of data are erroneous or uncertain, i.e. some data of each
 population do not deal with the population we want to compare. To give an example of such a situation with lack of information, we consider two
 populations -  New York and California people - reduced to people that take bus/trolley bus or walk to go working. Focusing on working people that
 take bus/trolley bus, suppose you are interested to know whether travel time of people from New York is significantly different to the one from 
California from a sample of people where the place they live -  New York or California - is available but the way of travel (the {\it{label}}) associated with each data in hand is not. \\
This kind of situation is the one we are interested in.  Indeed, we want to address the problem of testing means of (sub)populations when the labels of data  are uncertain. More precisely, we first propose to show that this testing problem can be reformulated as testing the expected values of components from two samples of independent mixture variables. In our study of real data, we shall assume that the mixing-weights are known. It means that proportions of people walking or using  bus/trolley bus for each population (New York and California) are known, with respect to an auxiliary variable (age for instance).  Then, we provide a testing procedure that takes into account this  information on populations - and we discuss about its performances. \\
The testing procedure we propose is directly inspired from ideas in Autin and Pouet \cite{AuPo}. In this previous work, a nonparametric procedure has been proposed to test whether the densities of two independent samples of independent random variables result from the same mixture of components or not. The value of the test statistic requires to invert in some sense the mixing-weights operators of samples (see Definition \ref{def-op})  as a preliminary step  to be  calculated. This testing procedure was proved to be powerful since it is minimax over Besov spaces (more details are given in paragraph $3.1$ in Autin and Pouet \cite{AuPo}). 
More focusing on practical purposes, we show that providing a testing procedure that incorporates combinatory ideas -  provided that the mixing-weights are known - is quite relevant compared to a procedure  usually used in classification.\\

\noindent Paper is organized as follows. In Section \ref{two} we present the mixture-model we are interested in. Connection between the testing problem for which the labels of data are not certain and the problem of testing the expected value of the components  involved in the mixture-model is provided.  In Section \ref{three}, we present three testing procedures. The first one is the {\it{Oracle Procedure}} that uses Welch's t-test on data associated with the label of the components we want to test. Of course this procedure is not tractable for the testing problem with lack of information about labels but it will be used as a benchmark to assess the loss of performances of the other tractable testing procedures. The second testing procedure  we present is the {\it{Expert Procedure}} that uses Welch's t-test on data that are supposed to have, with probability larger than or equal to one half,  the label associated with the components we focus on. The third and last procedure, namely the {\it{Mixing Procedure}}, uses  combinatory properties leading to a new performing test.  Section \ref{four} deals with numerical experiments to point out  the good performances of the Mixing Procedure - the one we suggest - compared to the Expert Procedure and to assess the loss of performances due to the mixing-effect compared to the Oracle Procedure. An application to real data is also presented whereas a brief conclusion is postponed in Section \ref{six}. Finally, the technical lemmas and the proposition we used to prove our main theoretical result (see Theorem \ref{main}) together with their proofs  can be found in the appendix.

\section{Model description  and hypothesis testing problem}\label{two}

\subsection{Mixture-models with varying mixing-weights}
\noindent Let $X_{1},\ldots, X_{n}$ be independent random variables such that, for any $1 \leq i \leq n,$ the density of $X_{i}$ on $\mathbb{R}$, denoted by $f_{_{X_{i}}}$, is a
 mixture density with 
components $p_1$ and $p_2$ and mixing-weights $\omega_1(i)$ and 
$\omega_2(i)$, i.e. 
$$ f_{_{X_{i}}}=\omega_1(i) p_1 + \omega_2(i) p_2.$$
We also introduce labels attached to $X_{1},\ldots,X_{n}$, denoted by
$u_{1},\ldots,u_{n}$.  
This point of view is one interpretation of mixture-models among others (see
Section 1.4 in McLachlan and Peel \cite{McLaPe00}). The main difference lies in considering
varying mixing-weights in our model. This point is very important  (see Autin 
and Pouet\cite{AuPo}). Therefore our model cannot be described as a mixture-model in
the usual sense.\\

\noindent Similarly to the sample $X_{1},\ldots, X_{n}$, we consider a sample
of independent random variables $Y_{1},\ldots, Y_{n'}$ such that, for any $1 \leq i \leq n$ the density of $Y_{i}$ on $\mathbb{R}$, denoted by $f_{_{Y_{i}}}$, is a
 mixture density with 
components $p'_1$ and $p'_2$ and mixing-weights $\omega'_1(i)$ and 
$\omega'_2(i)$, i.e. 
$$ f_{_{Y_{i}}}=\omega'_1(i) p'_1 + \omega'_2(i) p'_2.$$
We also introduce labels attached to $Y_{1},\ldots,Y_{n}$, denoted by
$v_{1},\ldots,v_{n}$ and we assume that this second sample is independent from the first one.\\

\noindent If \ $^t.$ denotes the transpose operator, the two mixture-models we have just introduced can be rewritten in a simpler way as follows:
\begin{eqnarray}\label{mix-mod}
{\bf{f_{_{X}}}}=\Omega_{_{X}} {\bf{p}} \quad \hbox{ and } \quad {\bf{f_{_{Y}}}}=\Omega_{_{Y}} {\bf{p'}},
\end{eqnarray} 

\noindent with,  
\begin{itemize}
\item[-] ${\bf{f_{_{X}}}}= {^t(f_{_{X_{1}}}, \dots, f_{_{X_{n}}})}$, \quad ${\bf{f_{_{Y}}}}= {^t(f_{_{Y_{1}}}, \dots, f_{_{Y_{n}}})}$,
\item[-] $ {\bf{p}} = {^t(p_{1},p_{2})}$, \quad $ {\bf{p'}} = {^t(p'_{1},p'_{2})}$,
\item[-] $\Omega_{_{X}}=\left(\omega_{l}(i)\right)_{i,l}$, \quad  $\Omega_{_{Y}}=\left(\omega'_{l}(i)\right)_{i,l}$.
\end{itemize}

\noindent \begin{definition}\label{def-op}
The matrices $\Omega_{_{X}}$ and $\Omega_{_{Y}}$ involved in the model (\ref{mix-mod}) are called {\it{the mixing-weights operators}}.
\end{definition}
\noindent \begin{definition}\label{def-var-mix}
Any mixture-model (\ref{mix-mod}) such that $\Omega_{_{X}}$ and $\Omega_{_{Y}}$ are full rank matrices is called mixture-model with varying mixing-weights.  
\end{definition}

\subsection{Example of modeling with mixture-models}
\noindent Let us illustrate this theoretical set-up with the example cited 
in the introduction. The random variables $X_{1},\ldots, X_{n}$ 
correspond to the travel times of people in the state of New York and
the random variables $Y_{1},\ldots, Y_{n}$ to
travel times in the state of California. The labels are the ways of transportation to go working and can be
either \textit{Bus/trolley bus} (label $1$) or \textit{Walk} (label $2$) . 
The last step to complete the mixture 
model is to describe the mixing-weights for each observation. 
In each state the mixing-weights strongly depend on the age (over $21$ or under $20$ years old (y.o.)). Table \ref{tab:label0} illustrates  this fact.

\begin{table}
  \caption{  \label{tab:label0}
Populations weights (and sizes) with respect to age}
  \centering
\begin{tabular}[h]{|c|c|c|} \hline 
 & Bus/trolleybus & Walk \\ \hline \hline
New York over 21 y.o.& 51.93\% (4313) & 48.07\% (3993)  \\ 
New York under 20 y.o.& 34.65\% (306)  & 65.35\% (577) \\ \hline
California over 21 y.o.& 57.4\% (4479) & 42.6\% (3324)  \\ 
California under 20 y.o.& 42.77\% (497) & 57.23\% (665)\\ \hline
\end{tabular}
\end{table}

\noindent This table leads to the following mixing-weights:
\begin{eqnarray*}
\omega_1(i) & = & 0.5193, \omega_2(i) = 0.4807 \\
 & & \mbox{ if the person $i$ is over 21 y.o. and   lives in New York}, \\
\omega_1(i) & = & 0.3465, \: \omega_2(i) = 0.6535 \\
 & & \mbox{ if the person $i$ is under  20 y.o. and   lives in New York}, \\
\omega'_1(i) & = & 0.574, \: \omega'_2(i) = 0.426 \\
 & & \mbox{ if the person $i$ is over 21 y.o. and  lives in California}, \\
\omega'_1(i) & = & 0.4277, \:  \omega'_2(i) = 0.5723  \\
 & & \mbox{ if the person $i$ is over 21 y.o. and   lives in California}.
\end{eqnarray*}

\noindent The reader can legitimatelly wonder why the age is assumed to be known and 
not the ways of transportation to go working. One can think about at least two good
 reasons. The first
one is an a priori reason. The survey would be rather lengthy if all
interesting variables were included.  Therefore the survey is restricted to 
a small set of 
informative variables strongly linked with the interesting variables. Moreover
these informative variables can be chosen as objective as possible and
easily recordable.
This can be called planned missing values (see Graham \cite{Graham}).
The other reason is an a posteriori one. During the data analysis of a survey, 
researchers are often confronted with new hypotheses to test. In many situations, the relevant variables have not been recorded and researchers have to plan 
a new survey which includes these new variables in order to check these hypotheses.
This leads to a waste of time and money.

\noindent Our testing problem on means from data with undefined labels can be associated with the testing problem (\ref{def-test-hyp}) in the mixture-model (\ref{mix-mod}). Indeed, it corresponds to a testing problem on means for which labels of data are unavailable: the only information on the $X_{i}$'s label (resp. $Y_{i}$'s label) is the probability $\omega_{l}(i)$ (resp. $\omega'_{l}(i)$) that
it corresponds to $l$, for any $l \in \{1,2\}$.   In other words, the added information on subpopulations we get is the knowledge of the mixing-weights operators.

\subsection{Hypothesis testing problem}
\noindent We recall that two data samples ${\bf{X}}= \ ^t(X_{1},\dots, X_{n})$ and ${\bf{Y}}= \ ^t(Y_{1},\dots, Y_{n})$ are considered. For a chosen label $l\in \{1,2\}$, we are interested in testing whether components $p_{l}$ and $p'_{l}$ have the same expected value or not. We want to address this problem in a general context that is: the parameters of variance $\sigma^2_{k}$ and $\sigma^{'2}_{k}$ of the components $p_{k}$ and $p'_{k}$ are unknown whatever $k \in \{1,2\}$.\\

\noindent For a fixed $l \in \{1,2\}$, when respectively denoting by $m_{l}$ and $m'_{l}$ the expected value of the components we focus on, the testing problem we consider is lying on the two following hypotheses: 
\begin{eqnarray}\label{def-test-hyp}
\hbox{the
 {\it{null hypothesis}} } &\mathcal{H}_{0}: \ m_{l}=m'_{l}, \nonumber \\
 &\\
 \hbox{the {\it{alternative hypothesis}} } &\mathcal{H}_{1}: \ m_{l}\not=m'_{l}. \nonumber
 \end{eqnarray}

\noindent We recall that providing a procedure to solve the testing problem
(\ref{def-test-hyp}) means giving a decision rule (or test) $\Delta \in \{0,1\}$ that relies on the value of a measurable function $T$  (test statistic) of $X_1,\dots, X_n$
and $Y_1,\dots, Y_n$.\\

\noindent  As usual, 
$\Delta=1$ will mean deciding $\mathcal{H}_{1}$ whereas $\Delta=0$ will mean deciding  $\mathcal{H}_{0}$.\\

\section{Description of  testing procedures}\label{three}

\noindent In this section we introduce the testing procedures we are interested in. 

\subsection{Oracle test: $\Delta_{o}$}

\noindent The first testing procedure we present is called the  {\it{Oracle Procedure}}. This is a two steps procedure. First step consists in recovering the true labels of data. Second step lies on using the Welch's t-test on data with label $l$ in order to know whether $m_{l}$  and $m'_l$ can be judged as different. This test  cannot be used in our context where the true labels are unknown but it will be used as a benchmark when comparing the performances of the other testing procedures. It corresponds to the procedure proposed by the {\it{oracle}}: any statistician having information on labels. \\

\noindent Here we describe into details the Oracle Procedure.
\noindent  Let us denote by 
\begin{itemize}
\item[-] $n_{{l}}=\begin{displaystyle}\sum_{i=1}^n\end{displaystyle}{\bf{1}}\{u_{i}=l\}$ and $n'_{{l}}=\begin{displaystyle}\sum_{i=1}^n\end{displaystyle}{\bf{1}}\{v_{i}=l\}$,
\item[-] $\bar{X}^{(l)}=\frac{1}{n_{{l}}}\begin{displaystyle}\sum_{i =1}^n\end{displaystyle}X_{i}{\bf{1}}\{u_{i}=l\}$ and $\bar{Y}^{(l)}=\frac{1}{ n'_{{l}}}\begin{displaystyle}\sum_{i=1}^n\end{displaystyle}Y_{i}{\bf{1}}\{v_{i}=l\}$,
\item[-] $\hat{\sigma}^2_{l}=\begin{displaystyle}\frac{1}{n_l}\sum_{i=1}^n\end{displaystyle}(X_{i}-\bar{X}^{(l)})^2{\bf{1}}\{u_{i}=l\}$ \hbox{ and } $\hat{\sigma}^{'2}_{l}=\begin{displaystyle}\frac{1}{n'_l}\sum_{i=1}^n\end{displaystyle}(Y_{i}-\bar{Y}^{(l)})^2{\bf{1}}\{v_{i}=l\}$.
\end{itemize}

\noindent The Oracle test $\Delta_{o}$ lies on the test statistic $T_{o}$ defined as follows
\begin{eqnarray*}
T_{o}&:=&\frac{|\bar{X}^{(l)}-\bar{Y}^{(l)}|}{\sqrt{\frac{\hat{\sigma}_{l}^2}{n_{l}} + \frac{\hat{\sigma}_{l}^{'2}}{n'_{l}}}}.
\end{eqnarray*}

\noindent Under the null hypothesis, the asymptotic law of $T_{o}$ is  known to be the Standard Gaussian one, namely $\mathcal{N}(0,1)$. Hence, 
$\Delta_{o}= {\bf{1}}\{T_{o}>q_{r}\}$ is a test with asymptotically  type $I$ error equal to $r$  $(0 < r <1)$, where $q_{r}$ is the quantile of order $1-\frac{r}{2}$ of the Standard Gaussian law.\\

\subsection{Expert test: $\Delta_{e}$}

\noindent The testing procedure we describe now is lying on a method used in classification. It is a two steps  procedure. The first step consists in allocating label $l$ to any data $X_{i}$ such that $\omega_{l}(i)\geq \frac{1}{2}$ and  to any data $Y_{j}$ such that $\omega'_{l}(j)\geq \frac{1}{2}$. The second step consists in using the Welch's t-test on the two subsamples of data that have been assigned to label $l$  to know whether $m_l$  and $m'_{l}$ can be judged as different. Notice that it means that the Welch's t-test is done on data having  possible wrong labels. \\

\noindent Put
\begin{itemize}
\item[-] $n_{_{l,e}}=\begin{displaystyle}\sum_{i=1}^n\end{displaystyle}{\bf{1}}\{\omega_{l}(i)\geq \frac{1}{2}\}$ and $n'_{_{l,e}}=\begin{displaystyle}\sum_{i=1}^n\end{displaystyle}{\bf{1}}\{\omega'_{l}(i)\geq \frac{1}{2}\}$,
\item[-] $\bar{X}^{(l)}_{e}=\frac{1}{n_{_{l,e}}}\begin{displaystyle}\sum_{i =1}^{n}\end{displaystyle}X_{i}.{\bf{1}}\left\{\omega_{l}(i)\geq \frac{1}{2}\right\}$ and $\bar{Y}^{(l)}_{e}=\frac{1}{n'_{_{l,e}}}\begin{displaystyle}\sum_{i =1}^{n}\end{displaystyle}Y_{i}.{\bf{1}}\left\{\omega'_{l}(i)\geq \frac{1}{2}\right\}$,
\item[-] ${\hat{\sigma}}_{_{l,e}}^2=\begin{displaystyle}\frac{1}{n_{_{l,e}}}\sum_{i =1}^{n}\end{displaystyle}(X_{i}-\bar{X}_e^{(l)})^2{\bf{1}}\left\{\omega_{l}(i)\geq \frac{1}{2}\right\}$,
\item[-] ${\hat{\sigma}}_{_{l,e}}^{'2}=\begin{displaystyle}\frac{1}{n'_{_{l,e}}}\sum_{i =1}^{n}\end{displaystyle}(Y_{i}-\bar{Y}_e^{(l)})^2{\bf{1}}\left\{\omega'_{l}(i)\geq \frac{1}{2}\right\}$.
\end{itemize}

\noindent The Expert test $\Delta_{e}$ relies on the test statistic $T_{e}$ defined as follows
\begin{eqnarray*}
T_{e}&:=&\frac{|\bar{X}_{e}^{(l)}-\bar{Y}_{e}^{(l)}|}{\sqrt{\frac{\hat{\sigma}_{_{l,e}}^2}{n_{_{l,e}}} + \frac{\hat{\sigma}_{_{l,e}}^{'2}}{n'_{_{l,e}}}}}.
\end{eqnarray*}

\noindent Then, the decision rule is done by putting
$\Delta_{e}= {\bf{1}}\{T_{e}>q_{r}\}$.\\

\subsection{Mixing test: $\Delta_{m}$}

\noindent The last testing procedure we propose is inspired from some ideas provided in Autin and Pouet \cite{AuPo}. Using combinatory methods, it proposes to invert in some sense the mixing-weights operators  so as to provide a new test that will be proved to perform well. Let us describe this new testing procedure into details.\\

\noindent Let us denote by $A_{_{X}}$ and  $A_{_{Y}}$ the matrices with $n$ lines and $2$ columns satisfying
\begin{eqnarray}
^t\Omega_{_{X}}A_{_{X}}=^t\Omega_{_{Y}}A_{_{Y}}=\begin{pmatrix}
n&0 \\
0&n
\end{pmatrix} \label{algebra}.\end{eqnarray}

\noindent {\underline{Notations}}: For any $(i,l) \in \{1,\dots,n\}\times \{1,2\}$, we denote respectively by $a_{l}(i)$ (resp. $a'_{l}(i)$) the entries of $A_{_{X}}$ (resp. $A_{_{Y}}$) associated with line $i$ and column $l$. \\

\noindent Following Maiboroda \cite{Mai-3} or Pokhyl'ko \cite{Pok-1}, solutions of equations (\ref{algebra}) are given by
\begin{eqnarray*}
a_{l}(i)&=& \frac{n}{det(^t\Omega_{_X}\Omega_{_X})}\sum_{k=1}^2(-1)^{l+k}\gamma^{(X)}_{_{lk}}\omega_k(i),\\
a'_{l}(i)&=& \frac{n}{det(^t\Omega_{_Y}\Omega_{_Y})}\sum_{k=1}^2(-1)^{l+k}\gamma^{(Y)}_{_{lk}}\omega'_{k}(i),\\
\end{eqnarray*}
\noindent where $\gamma^{(X)}_{_{lk}}$ and  $\gamma^{(Y)}_{_{lk}}$ are respectively the minor $(l,k)$ of the matrix $^t \Omega_{_X}\Omega_{_X}$ and  of the matrix $^t\Omega_{_Y}\Omega_{_Y}$. \\

\noindent For any $l \in \{1,2\}$, $\left(m_{l},m'_{l}\right)$ can be estimated  by the method of moments when using estimators $\left(\hat{m}_{l},\hat{m}'_{l}\right)$  defined as follows:
\begin{eqnarray*}
\left(\hat{m}_{l},\hat{m}'_{l}\right)&=&\left(\langle A_{_{X}}^{(l)},X \rangle_{n},\langle A_{_{Y}}^{(l)},Y \rangle_{n}\right)\\
&:=&\left(\frac{1}{n}\sum_{i=1}^na_{l}(i)X_{i},\frac{1}{n}\sum_{i=1}^n a'_{l}(i)Y_{i}\right),
\end{eqnarray*}

\noindent provided that $A_{_{X}}^{(l)}$ and $A_{_{Y}}^{(l)}$ respectively denote the $l$-th column-vector of matrices $A_{_{X}}$ and $A_{_{Y}}$.\\

\noindent The Mixing test $\Delta_{m}$ lies on the test statistic $T_{m}$ defined by:
\begin{eqnarray}\label{test-statistic}
T_{m}&:=&\frac{|\hat{m}_{l}-\hat{m}'_{l}|}{\sqrt{\hat{\mathbb{V}}_n^{(l)}}},
\end{eqnarray}
\noindent where $\hat{\mathbb{V}}_n^{(l)}$ is the estimated variance of $\hat{m}_{l}-\hat{m}'_{l}$, that is \\

\begin{eqnarray*}
\hat{\mathbb{V}}_n^{(l)}&=&\frac{1}{n^2} \sum_{i=1}^n \left[a_{l}^{2}(i)  \left( X_i - \omega_{1}(i) \hat{m}_l -  \omega_{2}(i)\hat{m}_l \right)^2 + {a'_{l}}^{2}(i)  \left( Y_i - \omega'_{1}(i) \hat{m}'_l -  \omega'_{2}(i)\hat{m}'_l \right)^2\right].\\
\end{eqnarray*}

\begin{rem} As discussed in Autin and Pouet \cite{AuPo}, for any $l \in \{1,2\}$ the random variable $\hat{m}_{l}$ (resp. $\hat{m}'_{l}$) is a good estimator for $m_{l}$ (resp. $m'_{l}$). Hence if the distance between $\hat{m}_{l}$  and $\hat{m}'_{l}$ is judged too large, the rejection of the null hypothesis $\mathcal{H}_{0}$ looks better. This idea motivates the choice of the test statistic $T_{m}$ we defined above. \end{rem}

\noindent Under the null hypothesis $\mathcal{H}_0$, the asymptotic law of $T_{m}$ is known, according to the following Theorem.

\begin{thm}\label{main}
\noindent Let $l \in \{1,2\}$. Assume that 
\begin{itemize}
\item the components within the mixture-model (\ref{mix-mod}) have moments with order $4$,
\item the mixing-weights of the mixture-model (\ref{mix-mod}) are such that 
\begin{eqnarray} \label{ass-sup}
\lim_{n \to +\infty} \frac {\begin{displaystyle}\sup_{i=1,\ldots,n} \end{displaystyle}a^2_l(i)} {\begin{displaystyle}\sum_{i=1}^n \end{displaystyle} a^2_l(i)} =\lim_{n \to +\infty} \frac {\begin{displaystyle}\sup_{i=1,\ldots,n} \end{displaystyle}{a}^{'2}_l(i)} {\begin{displaystyle}\sum_{i=1}^n \end{displaystyle} {a}^{'2}_l(i)}= 0.
\end{eqnarray}

\end{itemize}

\noindent Then, under the null hypothesis $\mathcal{H}_0$, the law of $T_{m}$ is asymptotically the Standard  Gaussian one, i.e.
\begin{eqnarray}\label{conv-loi-N}
T_{m} \stackrel{\mathcal{L}}{\to}\mathcal{N}(0,1).
\end{eqnarray} \end{thm}
 
 \noindent Hence, 
$\Delta_{m}= {\bf{1}}\{T_{m}>q_{r}\}$ is a test with asymptotically  type $I$ error equal to $r$ $(0 < r <1)$.

\begin{rem}
A wide range of mixing-weights of the mixture-model (\ref{mix-mod}) satisfy the condition (\ref{ass-sup}). Examples of such mixing-weights are given in  $(\ref{sure-ness})$ of Section \ref{four}.   \end{rem}

\noindent {\underline{Proof}}:\\

\noindent To prove Theorem \ref{main}, notice that it suffices to prove that for any $l \in \{1,2\}$
\begin{enumerate}
\item  $\frac{\begin{displaystyle}\frac{1}{n}\sum_{i=1}^n \end{displaystyle}a_l(i)X_i-m_{l}}{\sqrt{ \begin{displaystyle}\frac{1}{n^2}\sum_{i=1}^n\end{displaystyle} a_{l}^{2}(i)  \left( X_i - \omega_{1}(i) \hat{m}_1 -  \omega_{2}(i)\hat{m}_2 \right)^2}} \stackrel{\mathcal{L}}{\to}  \mathcal{N}(0,1)$,
\item  $\frac{\begin{displaystyle}\frac{1}{n}\sum_{j=1}^n \end{displaystyle}a'_l(i)Y_j-m'_{l}}{\sqrt{ \begin{displaystyle}\frac{1}{n^2}\sum_{j=1}^n\end{displaystyle} {a'_{l}}^{2}(i)  \left( Y_j - \omega'_{1}(j) \hat{m}'_1 -  \omega'_{2}(i)\hat{m}'_2 \right)^2}} \stackrel{\mathcal{L}}{\to}  \mathcal{N}(0,1)$.
\end{enumerate}

\noindent Because of the independence between the two samples and the fact that under the null hypothesis $\mathcal{H}_{0}$, $m_{l}=m'_{l}$. Since these two results of convergence can be proved by an analogous way, we only focus on proving the first one that can be rewritten as follows for any $l \in \{1,2\}$:
$$\frac{\begin{displaystyle}\sum_{i=1}^n \end{displaystyle}a_l(i)\left(X_i-\omega_1(i)m_1-\omega_2(i)m_2\right)}{\sqrt{ \begin{displaystyle}\sum_{i=1}^n\end{displaystyle} a_{l}^{2}(i)  \left( X_i - \omega_{1}(i) \hat{m}_1 -  \omega_{2}(i)\hat{m}_2 \right)^2}} \stackrel{\mathcal{L}}{\to}  \mathcal{N}(0,1).$$

\noindent  Denote, for any $n \in \mathbb{N}^*$, any $l \in \{1,2\}$ and any $1 \leq i \leq n$ 
\begin{eqnarray}
B_n^{(l)}&=& \sum_{i=1}^n a_{l}^{2}(i)  \mathbb E\left[ \left( X_i - \omega_{1}(i) m_1 -  \omega_{2}(i)m_2 \right)^2 \right],\label{B_n}\\
\hat{B}_n^{(l)}&=& \sum_{i=1}^n a_{l}^{2}(i)  \left( X_i - \omega_{1}(i) \hat{m}_1 -  \omega_{2}(i)\hat{m}_2 \right)^2. \label{hat_B_n}
\end{eqnarray} 

\noindent From Proposition \ref{prop-N},
\begin{eqnarray*}
\frac{\begin{displaystyle}\sum_{i=1}^n \end{displaystyle}a_l(i)\left(X_i-\omega_1(i)m_1-\omega_2(i)m_2\right)}{\sqrt{B^{(l)}_n}} \stackrel{\mathcal{L}}{\to}  \mathcal{N}(0,1).
\end{eqnarray*} 

\noindent In the sequel, we aim at proving that a same kind of result holds when replacing parameter $B^{(l)}_n$ by the estimator $\hat{B}^{(l)}_n$. In other words, the result we want to prove is the following:
\begin{eqnarray}\label{N-cons}
\frac{\begin{displaystyle}\sum_{i=1}^n \end{displaystyle}a_l(i)\left(X_i-\omega_1(i)m_1-\omega_2(i)m_2\right)}{\sqrt{\hat{B}^{(l)}_n}} \stackrel{\mathcal{L}}{\to}  \mathcal{N}(0,1).
\end{eqnarray}

\noindent  From Slutsky theorem, it  suffices to prove that estimator $\hat{B}^{(l)}_n$ of $B^{(l)}_n$ is consistent. We propose to divide the proof of this consistency
into two steps. First we prove
\begin{eqnarray}\label{conv-step1} \frac {\begin{displaystyle}\sum_{i=1}^n \end{displaystyle}a^2_l(i) \left( X_i - \omega_1(i) m_1 - \omega_2(i) m_2 \right)^2} {B^{(l)}_n}
\stackrel{Proba}{\to} 1.
\end{eqnarray}
The second step consists in replacing $m_1$ and $m_2$ by their consistent estimators $\hat{m}_1$ and $\hat{m}_2$ and in checking that the convergence in probability  still holds.\\

\noindent From this point we need more assumptions, that is to say the existence of the fourth order moment for $p_1$ and $p_2$.\\

\noindent Let us prove the first step. We apply Bienayme-Chebyshev inequality, for any $\epsilon>0$:
\begin{eqnarray*}
& & \mathbb P\left( \left| \begin{displaystyle}\sum_{i=1}^n \end{displaystyle}a^2_l(i) \left( X_i - \omega_1(i) m_1 - \omega_2(i) m_2 \right)^2-B_n^{(l)} \right| > B_n^{(l)}  \varepsilon\right) \\
&   \leq &\frac {\begin{displaystyle}\sum_{i=1}^n\end{displaystyle} a^4_l(i) \mathbb{V}ar\left(  \left( X_i - \omega_1(i) m_1 - \omega_2(i) m_2 \right)^2 \right)}
{(B_n^{(l)}\varepsilon)^2 }\\
&   \leq &\frac {\begin{displaystyle}\sum_{i=1}^n\end{displaystyle} a^4_l(i) \mathbb{E}\left[  \left( X_i - \mathbb{E}(X_{i}) \right)^4 \right]}
{(B_n^{(l)} \varepsilon)^2}\\
&   \leq &\frac{\begin{displaystyle}\hspace{-0.3cm}\sup_{\: \: j=1,\dots,n}\end{displaystyle}a^2_l(j)}{B_n^{(l)}}\: \frac{\begin{displaystyle}\sum_{i=1}^n\end{displaystyle} a^2_l(i) C(m_1,m_2,p_1,p_2)}
{B_n^{(l)} \varepsilon^2}\\
&   \leq &\frac{\begin{displaystyle}\hspace{-0.3cm}\sup_{\: \: j=1,\dots,n}\end{displaystyle}a^2_l(j)}{B_n^{(l)}}\: \frac{\left(\min(\sigma_1^2,\sigma_2^2)\right)^{-1}C(m_1,m_2,p_1,p_2)}
{\varepsilon^2}.
\end{eqnarray*}
\noindent Last inequalities  are obtained  by using Lemma \ref{lem-4} and Lemma \ref{lem-2}. The right part of the last inequality is the product of two terms. The left one tends to $0$ when $n$ goes to infinity because of assumption (\ref{ass-sup}). The right one is a constant that only depends on $\epsilon$ and the parameters of $p_1$ and $p_2$. When considering the limit in infinity with respect to $n$, we conclude that property (\ref{conv-step1}) holds.\\

\noindent We end by proving the second step. We have
\begin{eqnarray*}
& &\sum_{i=1}^n a^2_l(i) \left( X_i - \omega_1(i) \hat m_1 - \omega_2(i) \hat m_2 \right)^2 \\
& = & \sum_{i=1}^n a^2_l(i) \left( X_i - \omega_1(i)  m_1 - \omega_2(i)  m_2 \right)^2 \\
 & & \hspace{-0.4cm}+ 2 \sum_{i=1}^n a^2_l(i) \hspace{-0.1cm}\left( X_i \hspace{-0.1cm}- \hspace{-0.1cm}\omega_1(i) m_1 \hspace{-0.1cm}- \hspace{-0.1cm}\omega_2(i) m_2 \right)\hspace{-0.1cm} \left( \omega_1(i) (m_1 \hspace{-0.1cm}- \hspace{-0.1cm}\hat m_1) + \omega_2(i)
(m_2 \hspace{-0.1cm}- \hspace{-0.1cm}\hat m_2) \right) \\
& & \hspace{-0.3cm}+ \sum_{i=1}^n a^2_l(i) \left( \omega_1(i) (m_1 - \hat m_1) + \omega_2(i) (m_2 - \hat m_2) \right)^2.
\end{eqnarray*}
The first term is exactly the one appearing in the first step and also converges to $1$  in probability when divided by $B^{(l)}_n$.  We turn to the second term. Cauchy-Schwarz inequality entails that
\begin{eqnarray*}
& &  \hspace{-0.05cm} \left|\sum_{i=1}^n a^2_l(i) \left( X_i  \hspace{-0.1cm}-  \hspace{-0.1cm}\omega_1(i) m_1  \hspace{-0.1cm}-  \hspace{-0.1cm}\omega_2(i) m_2 \right) \left( \omega_1(i) (m_1  \hspace{-0.1cm}- \hspace{-0.1cm} \hat m_1) + \omega_2(i)
(m_2  \hspace{-0.1cm}- \hspace{-0.1cm} \hat m_2) \right)\right| \\
& \leq & \sqrt{  \sum_{i=1}^n a^2_l(i) \left( X_i - \omega_1(i) m_1 - \omega_2(i) m_2 \right) ^2} \\
& & \times
\sqrt{2 (m_1- \hat m_1)^2 \sum_{i=1}^n a^2_l(i)  \omega^2_1(i) +  2 (m_2 - \hat m_2)^2
 \sum_{i=1}^n a^2_l(i)  \omega^2_2(i)}.
\end{eqnarray*} 
When divided by $B^{(l)}_n$, the first term of the righthand-side converges to $1$ in probability: it is the result of the first step.  By using Lemma \ref{lem-2}, one gets 
$$ \max\left(\begin{displaystyle}\sum_{i=1}^n \end{displaystyle}a^2_l(i) \omega^2_1(i),
\begin{displaystyle}\sum_{i=1}^n \end{displaystyle}a^2_l(i) \omega^2_2(i)\right)\leq B^{(l)}_n\left(\min(\sigma_1^2,\sigma_2^2)\right)^{-1}.$$
\noindent Hence the second term of the right-hand side of the inequality converges to $0$ in probability when divided by $B^{(l)}_n$ because of the consistency of estimators $\hat{m}_l$ (see Lemma \ref{lem-5}).  Hence second the term we are interested in converges to $0$ in probability when divided by $B^{(l)}_n$. \noindent We can proceed in the same way in order to prove that the third term converges to $0$ in probability when divided by $B^{(l)}_n$.\\

\noindent So, we have just proved that 
\begin{eqnarray}\label{conv-step2} \frac {\hat{B}^{(l)}_n} {B^{(l)}_n}
\stackrel{Proba}{\to} 1.
\end{eqnarray}

\noindent We conclude that the exact variance $B^{(l)}_n$ can be replaced by the consistent estimator $\hat{B}^{(l)}_n$ for the result of convergence. In other words,  the property (\ref{N-cons}) holds.\\

\section{Numerical experiments}\label{four}

\subsection{Numerical performances of the Mixing-test}\label{four-1}
\noindent In this section we provide numerical experiments and we discuss about the performances of our testing procedure. What we often expect
is a gain of performance of the test $\Delta_m$ - that is to say a smaller  type $II$ error when  the type $I$ error is chosen to be $r=0.05$ -  comparatively to the test $\Delta_e$. Without loss of generality, we suppose that $n$ is even.\\
We consider the Gaussian setting and we assume in this section that the mixing-weights operators $\Omega_{_{X}}$ and $\Omega_{_{Y}}$ have the following form:
\begin{eqnarray}\label{sure-ness}\Omega_{_{X}}=\begin{pmatrix} \alpha & 1-\alpha \\\dots & \dots  \\ \alpha & 1-\alpha \\ 1-\beta & \beta \\ \dots & \dots \\  1-\beta &\beta  \end{pmatrix} \hbox{ and } \Omega_{_{Y}}=\begin{pmatrix} \alpha'& 1-\alpha' \\ \dots & \dots  \\ \alpha' & 1-\alpha' \\ 1-\beta'& \beta' \\ \dots & \dots \\  1-\beta' &\beta' \end{pmatrix},  \end{eqnarray}
where $\frac{n}{2}$ data  from ${\bf{X}}$ (resp. ${\bf{Y}}$) deal with the couple of mixing-weights $(\alpha,1-\alpha)$ (resp. $(\alpha',1-\alpha')$) and 
the other $\frac{n}{2}$ data from ${\bf{X}}$ (resp. ${\bf{Y}}$) deal with  the couple of mixing-weights $(1-\beta,\beta)$ (resp. $(1-\beta',\beta')$).
Suppose now that our testing problem is dealing with the first component, i.e. $l=1$ and that $\Omega_{_{X}}$and $\Omega_{_{Y}}$ are full rank matrices, i.e. $\alpha+\beta \not= 1$ and  $\alpha'+\beta' \not= 1$. 

\subsubsection{Mixing-test versus Expert-test}

\noindent In this paragraph we provide a motivation for the use our testing procedure $\Delta_{m}$. For the sake of simplicity we suppose that $\alpha=\beta$ and that $\alpha'=\beta'$. For any value of $(\alpha,\alpha') \in ]\frac{1}{2},1[^2$, there are many situations where the performance of the Expert test is quite bad even if the numbers of observations $n$ is large. 

\begin{itemize}
\item Dealing with two components with equal expected value, $\Delta_{e}$ can most of time detect a difference between these components (wrong decision) whereas our test doesn't. For instance, suppose that $m_{1}=m'_{1}$ and that $m'_{2}$ is large away  from $m_{2}$ as $\alpha=\alpha'$.  \noindent  Since 
$\mathbb{E}\left(\bar{X}^{(1)}_{e}\right)\not=\mathbb{E}\left(\bar{Y}^{(1)}_{e}\right),$ using $\Delta_{e}$ to detect equality between components $m_{1}$ and $m'_{1}$  would be a very bad choice in that context. For $n$ large enough, it would imply that $T_{e} > t_{r}$ with high probability.  Hence, the wrong decision $\mathcal{H}_1$ may often be done.\\

\noindent An example of such a situation is given here in the case where $\alpha=\alpha'=0.9$. Consider the testing problem (\ref{def-test-hyp}) and suppose $\sigma_{1}=\sigma'_{1}=\sigma_{2}=\sigma'_{2}=1$ and that $m_{1}=m'_{1}=0$,  $m_{2}=1$ and $m'_{2}=m_{2}+\delta $.  For varying values of $n$, $\delta$ and $40 \ 000$ repetitions of  $\Delta_{e}$ with $r=0.05$, we give the percentage of wrong decisions $\mathcal{H}_1$ in Table  \ref{tab:label1}. 

\begin{table}
  \caption{  \label{tab:label1}
Percentage of wrong decisions by $\Delta_{e}$}
  \centering
  \begin{tabular}[h]{|c|c|c|c|c|c|}
\hline 
    $\ \delta \ / \: \ n$ &  $100$ & $200$ & $500$ & $1000$ & $2000$   \\ 
%      \hline 
%$1.3$ & $0.0554$ & $0.0556$ & $0.0625$ & $0.0732$ & $0.0976$\\
  \hline 
$0.5$ & $0.057$ & $0.064$ & $0.086$ & $0.121$ & $0.191$\\
\hline 
$1$ & $0.074$ & $0.098$ & $0.172$ & $0.302$ & ${\textcolor{red}{0.521}}$\\
\hline 
$2$ & $0.126$ & $0.210$ & $0.462$ & ${\textcolor{red}{0.749}}$ & ${\textcolor{red}{0.963}}$\\
\hline 
$3$ & $0.188$ & $0.350$ & ${\textcolor{red}{0.722}}$ & ${\textcolor{red}{0.950}}$ & ${\textcolor{red}{0.999}}$\\
\hline
  \end{tabular}
\end{table}

\noindent Notice that the percentage of wrong decisions by $\Delta_{e}$ turns up as $n$ grows up and can be quite important if $m'_{2}$ is sufficiently far away from $m_{2}$. Most of time, the expert detects a difference between the components $m_{1}$ and $m'_{1}$ but there is not in that context. Comparatively speaking, the percentage of wrong decisions by $\Delta_{m}$ is around $0.05$.

\item Most of time $\Delta_{e}$  fails to detect a difference between two components with different expected value whereas our test doesn't. For instance, suppose that $m_{1}\not=m'_{1}$ and that $$m_{2}\approx (1-\alpha)^{-1}\left(\alpha' m'_{1}+(1-\alpha') m'_{2}-\alpha m_{1}\right).$$ \noindent  Since 
$$\mathbb{E}(X_{i}) \approx \mathbb{E}(Y_{i}), \quad \hbox{ for any } 1 \leq i\leq \frac{n}{2},$$ using $\Delta_{e}$ to detect the difference between  $m_1$ and $m'_1$  would be a very bad choice in that context. Indeed, according to the law of large numbers, with high probability  - that increases as $n$ goes up - $\bar{X}_{e}^{(1)}$ and $\bar{Y}^{(1)}_{e}$ would be very close to each other.  It would imply that $T_{e} \leq 1.96$ with high probability.  True decision $\mathcal{H}_1$ would be taken only in $5 \ \%$ of cases.\\

\noindent An example of such a situation is given here in the case $\alpha=\alpha'=0.9$. Consider the testing problem (\ref{def-test-hyp}) and suppose that $\sigma_{1}=\sigma'_{1}=\sigma_{2}=\sigma'_{2}=1$ and that $m_{1}=0$, $m'_{1}=0.1,$ $m_{2}=1$ and $m'_{2}=2$.  For varying values of $n$ and $40 \ 000$ repetitions of  $\Delta_{m}$ with $r= 0.05$, we give the percentage of correct decisions from $\Delta_m$ in Table  \ref{tab:label2}. \\

\begin{table}
 \caption{ \label{tab:label2} Percentage of correct decisions by $\Delta_{m}$}
  \centering
  \begin{tabular}[h]{|c|c|c|c|c|c|c|c|}
\hline 
    $n$ &  $500$ & $1000$ & $2000$ & $3000$ & $4000$ & $5000$ & $6000$  \\ 
\hline 
 $\Delta_{m}$ & $0.146$ & $0.242$ & $0.438$ & ${\textcolor{red}{0.595}}$ & ${\textcolor{red}{0.718}}$ & ${\textcolor{red}{0.810}}$ & ${\textcolor{red}{0.879}}$\\
\hline 
  \end{tabular}
 \end{table}

\noindent As expected, the percentage of correct decisions by $\Delta_{m}$ goes up as $n$ grows up. But it is not the case for the percentage of correct decisions by $\Delta_{e}$ which are always around $0.05$. Most of time, the expert is unable to detect the difference between the components in that context. 
\end{itemize}

\noindent Finally we conclude that it is better to choose $\Delta_{m}$ for the problem we are interested in. 

\subsubsection{Mixing-test versus Oracle-test}

\noindent In this paragraph we compare the {\it{empirical powers}}  of $\Delta_{m}$  to the Oracle test $\Delta_{o}$ ones, when taking $r=0.05$ and the same parameters as the last example.
We recall that the empirical power  of any test $\Delta$ corresponds to the numerical evaluation of the probability to correctly decide $\mathcal{H}_1$, according to $\Delta$.\\

\begin{table}
\caption{ \label{tab:label3} Empirical Power of  $\Delta_{o}$ and $\Delta_{m}$}
 \centering
  \begin{tabular}[h]{|c|c|c|c|c|c|c|c|}
\hline 
   Test \ / \: \ $n$ &  $500$ & $1000$ & $2000$ & $3000$ & $4000$   & $5000$ & $6000$\\ 
 \hline 
$\Delta_{o}$ & $0.200$ & $0.349$ & $0.609$ & $0.783$ & $0.886$ & $0.942$ & $0.973$\\
\hline
$\Delta_{m}$ & $0.149$ & $0.245$ & $0.427$ & ${\textcolor{red}{0.585}}$ & ${\textcolor{red}{0.704}}$ & ${\textcolor{red}{0.798}}$ & ${\textcolor{red}{0.868}}$\\
\hline
  \end{tabular}
\end{table}

\noindent According to Table  \ref{tab:label3}  we remark that  the bigger $n$ the better the powers of  $\Delta_{o}$ and of $\Delta_{m}$.  Moreover we note that the empirical power of the Mixing test is not bad when comparing to the Oracle test. \\

\vspace{0.5cm}

\noindent In Table  \ref{tab:label4} we give the empirical power $\mathcal{P}_{m}$ of $\Delta_m$ measured in samplings of size $n=1000$ in the case where $m_{1}=0$, $m'_{1}=\bar{\delta}$, $m_{2}=1$ and $m'_{2}=0$. \\

\begin{table}
\caption{ \label{tab:label4} Empirical Power $\mathcal{P}_{m}$ of  $\Delta_{m}$ for varying $\alpha=\alpha'$}
 \centering
  \begin{tabular}[h]{|c|c|c|c|c|c|c|c|c|c|}
\hline 
 $\bar{\delta}$ \ / \: \ $\alpha$ &  $0.6$ & $0.65$ & $0.7$ & $0.75$ & $0.8$   & $0.85$ & $0.9$ & $0.95$ & $1$\\ 
 \hline 
  $0.1$ &  $0.071$ & $0.092$ & $0.125$ & $0.161$ & $0.196$   & $0.239$ & $0.275$ & $0.317$ & $0.353$\\ 
   \hline 
  $0.2$ &  $0.130$ & $0.226$ & $0.353$ & $0.484$ & ${\textcolor{red}{0.603}}$   & ${\textcolor{red}{0.700}}$ & ${\textcolor{red}{0.780}}$ & ${\textcolor{red}{0.837}}$ & ${\textcolor{red}{0.883}}$\\ 
   \hline 
  $0.3$ &  $0.238$ & $0.446$ & ${\textcolor{red}{0.659}}$ & ${\textcolor{red}{0.816}}$ & ${\textcolor{red}{0.912}}$   & ${\textcolor{red}{0.960}}$ & ${\textcolor{red}{0.983}}$ & ${\textcolor{red}{0.999}}$ & ${\textcolor{red}{1}}$\\ 
   \hline 
  $0.4$ &  $0.374$ & ${\textcolor{red}{0.671}}$ & ${\textcolor{red}{0.882}}$ & ${\textcolor{red}{0.938}}$ & ${\textcolor{red}{0.993}}$   & ${\textcolor{red}{1}}$ & ${\textcolor{red}{1}}$ & ${\textcolor{red}{1}}$ & ${\textcolor{red}{1}}$\\ 
  \hline
  \end{tabular}
\end{table}

\noindent As expected, looking at  Table \ref{tab:label4}, 
\begin{itemize}
\item quantity $\mathcal{P}_{m}$ depends on the intrinsic difficulty of the problem. Indeed the larger the absolute value of quantity $\bar{\delta}:=m'_{1}-m_{1}$, the easier the problem of detection and so the more powerful the test,
\item the larger the degree of certainty $\alpha$ the better the power of $\Delta_{m}$. This is due to the fact that the expectation of the number of wrong labels considered by the Expert Procedure grows up as  $\alpha$ goes down.
\end{itemize}

\subsubsection{Comparisons on performances of $\Delta_{m}$ for varying values of $(\alpha,\alpha')$}

\noindent As previously discussed, we expect that the better the degree of certainties of the expert, the better the performance of the test $\Delta_{m}$. This statement is highlighted here when considering the same parameters as before, $\bar{\delta}=0.5$  and many choices of couple ($\alpha, \alpha'$).
For each choice of ($\alpha, \alpha'$) done, we provide in Table  \ref{tab:label5} the empirical power $\mathcal{P}_{m}$ of our test $\Delta_{m}$ that is the percentage of correct detection of a difference between $m_{1}$ and $m'_{1}$. 

\begin{table}
\caption{ \label{tab:label5} Percentage of correct decisions by $\Delta_{m}$}
 \centering
  \begin{tabular}[h]{|c|c|c|c|}
\hline 
    $(\alpha,\alpha')$,  $/$  $n$ &  $100$ & $200$ & $500$ \\
    \hline 
$(0.90,0.60)$ & $0.153$ & $0.254$ & $0.533$ \\ 
    \hline 
$(0.80,0.70)$ & $0.311$ & $0.536$ & $0.896$ \\ 
    \hline 
$(0.75,0.75)$ & $0.332$ & $0.564$ & $0.918$ \\ 
\hline 
  \end{tabular}
  \end{table}

\noindent Interpretation of the results presented in  Table \ref{tab:label5} goes in the same way as Autin and Pouet \cite{AuPo}: the  bigger the smallest eigenvalue of both operators  $^t \Omega_{_{X}}\Omega_{_{X}}$ and $^t \Omega_{_{Y}}\Omega_{_{Y}}$ that is $\lambda_{_{min}}=\frac{1}{2}\left(1-2\min(\alpha,\alpha')(1-\min(\alpha,\alpha'))\right)$  the better the power of our test $\Delta_{m}.$ Note that, the larger the minimum value between $\alpha \in ]\frac{1}{2},1[$ and $\alpha' \in ]\frac{1}{2},1[$ the bigger $\lambda_{_{min}}$. 

\subsubsection{Brief conclusion}

\noindent Let us summarize the main facts. First, in some cases experts can be completely wrong because of the overall design, that is to say the link between the means of the components and the mixing-weights. This is a serious issue for the Expert test. The results become worse and worse as the sample size increases. The test adapted to the varying mixing-weights that we propose does not suffer from this drawback. The second fact is the good behavior of our test compared to the Oracle test. Although it is behind, the power is quite acceptable. The last important fact which has already been stressed by Autin and Pouet \cite{AuPo} is the effect of the mixing-weights. It is known a priori thanks to the smallest eigenvalue of the operators $^t\Omega_{_{X}}\Omega_{_{X}}$ and $^t\Omega_{_{Y}}\Omega_{_{Y}}$. This point is important as the statistician can act in order to counter to this effect, e.g. he can improve the accuracy of the expert system giving the mixing-weights or increase the sample sizes.

\subsection{Application to real data}\label{five}

\noindent In this section we apply our methodology to real data and we
discuss about the results.
\subsubsection{Description of the data}
\noindent We have selected data from U.S. Census Bureau website, more precisely PUMS
 2006 (see \cite{US}). 
We are interested in comparing travel time of people living either
 in the state of New York (abbreviated in NY)
 or either in the state of California (abbreviated in CA). Two ways of
 transportation have been kept: Bus/trolley bus and Walk.
 We have also kept a variable linked to age as it will be useful for the 
 mixture-model with varying mixing weights. This variable records the fact that a person is
over $21$ years old or under $20$ years old.\\

\noindent Here are few facts to roughly describe the PUMS sample. Table \ref{tab:label5bis}
gives one-level information.

\begin{table}
\caption{\label{tab:label5bis} Description of the population}
\centering
\begin{tabular}[h]{|c|c|c|} \hline 
 & NY & CA \\ \hline \hline
Total & 9189 & 8935 \\ \hline
Over 21 & 90.39\% (8306) & 87.04\% (7803)  \\
Under 20 &  9.61\% (883) & 12.96\% (1162) \\ \hline
Walk &  49.73\% (4570) &   44.5\% (3989)\\ 
Bus/trolley bus & 50.27\% (4619)  & 55.5\% (4979) \\ \hline
\end{tabular}
\end{table}

\noindent In Table \ref{tab:label6} we compute the mean and the standard deviation (in parentheses) of  the travel time according to the categorical variable means of transportation to go working.

\begin{table}
\caption{\label{tab:label6} One-way analysis of travel time (in minutes)}
\centering
\begin{tabular}[h]{|c|c|c|c|} \hline 
 & Walk  & Bus/trolley bus & Walk and Bus/trolley bus \\ \hline \hline
NY &  12.25 (12.18)  & 47.26 (28.79) & 29.85  (28.23)\\
CA & 11.23  (12.23)  & 45.12 (28.84) &   30.04 (28.49) \\ \hline
\end{tabular}
\end{table}

\noindent As it can be seen in Table \ref{tab:label6} there might be no difference between
New York and California. Nevertheless if the means of transportation is unavailable,
it will be perilous to decide when considering the whole sample without any
other information. Indeed as shown in Table \ref{tab:label5bis}, the difference between 
New York and California is decreased because of the structure of the population (less
people under 20 years old in New York).

\subsubsection{Methodology}

\noindent We assume in the sequel that the information about the way of transport (labels) are unavailable at the microdata level.
We are going to apply the test $\Delta_m$ adapted to the varying mixing-weights mixture-model. The age variable is the only auxiliary information available at the microdata level that permits to get the mixing-weights to our mixture-model (\ref{mix-mod}). \\

\noindent For comparison purpose we have also applied the so-called Expert test. 
The type I error is chosen to be $0.1$. \\

\noindent According to the notations we introduced in (\ref{sure-ness}) and to Table \ref{tab:label0}, 
\begin{eqnarray}
(\alpha,\beta)=(0.5193,0.6535) & (\alpha',\beta')=(0.574,0.5723).
\end{eqnarray}

\noindent We consider the following sample:
 $500$ persons over $21$ and $500$ persons under $20$ were randomly sampled in
 each state ($n=1000$).\\

\noindent We applied three testing procedures:
\begin{enumerate}
\item Oracle test,
\item Expert test,
\item Mixing test.
\end{enumerate}

\noindent First we test the equality of the averages when the ways of
transportation to work is Bus/trolley bus (label $1$) in Table \ref{tab:label7}.
 In this case, the other means of
transportation to work is considered as a nuisance parameter.

\begin{table}
\caption{\label{tab:label7} Bus/trolleybus decisions}
\centering
\begin{tabular}[h]{|c|c|c|} \hline 
 & Decision $n=1000$ & p-value \\ \hline \hline
Oracle test & not rejected & $0.24$ \\
Expert test & not rejected & $0.42$ \\ \hline
Mixing test & not rejected & $0.11$ \\ 
\hline
\end{tabular}
\end{table}

\noindent Next we reverse the set-up. We test the equality of the 
averages when the means of 
transportation to work is Walk (label $2$) in Table \ref{tab:label8}. Bus/trolley bus is now a nuisance 
parameter.

\begin{table}
\caption{\label{tab:label8} Walk decisions}
\centering
\begin{tabular}[h]{|c|c|c|} \hline 
 & Decision $n=1000$  & p-value \\ \hline \hline
Oracle test & not rejected & $0.23$ \\
Expert test & rejected &  $0.04$ \\ \hline
Mixing test & not rejected & $0.48$    \\ 
\hline
\end{tabular}
\end{table}

\subsubsection{A tough situation}
\noindent Here we are also interested in comparing travel time of people living either
 in the state of New York or either in the state of Illinois (abbreviated in IL). Data come from U.S. Census Bureau \cite{US}. Two ways of
 transportation to work have been kept: Bus/trolley bus or Railroad.
 We have also kept the gender variable as it will be useful for the varying
 mixing-weights mixture-model. As it will be seen, the situation is much 
more involved compared to the one in the previous section.\\

\noindent Here are few facts to roughly describe the PUMS sample. Table \ref{tab:label9}
gives one-level information.

\begin{table}
\caption{\label{tab:label9} Description of the population}
\centering
\begin{tabular}[h]{|c|c|c|} \hline 
 & NY & IL \\ \hline \hline
Total & 6974 & 2899 \\ \hline
Men &  46.5\% (3247) & 48.2\% (1398)  \\
Women &  53.5\% (3727) & 51.8\% (1501) \\ \hline
Bus/trolley bus & 66.2 \% (4619)  & 58.4\% (1692) \\
Railroad &  33.8\% (2355) &   41.6\% (1207)\\ \hline
\end{tabular}
\end{table}

\noindent The mixing-weights depend on the gender as illustrated in Table \ref{tab:label10}.

\begin{table}
\caption{\label{tab:label10} Mixing-weights}
\centering
\begin{tabular}[h]{|c|c|c|} \hline 
 & Bus/trolley bus  & Railroad \\ \hline \hline
NY men & 55.8 \% (1813) & 44.2\% (1434)\\
NY women & 75.3 \% (2806) & 24.7\% (921)\\ \hline
IL men &  50.8 \% (710) & 49.2\% (688)  \\ 
IL women & 65.4 \% (982)& 34.6\% (519)     \\ 
\hline
\end{tabular}
\end{table}

\noindent According to the notations we introduced in (\ref{sure-ness}) and Table \ref{tab:label10}, 
\begin{eqnarray}
(\alpha,\beta)=(0.558,0.247) & (\alpha',\beta')=(0.508,0.346).
\end{eqnarray}

\noindent In Table \ref{tab:label11} we compute the mean and standard deviation (in parentheses) of the travel time according to the categorical variable way of transportation to work.

\begin{table}
\caption{\label{tab:label11} One-way analysis of travel time}
\centering
\begin{tabular}[h]{|c|c|c|c|} \hline 
 & Bus/trolley bus  & Railroad & Bus/trolley bus and Railroad \\ \hline \hline
New York & 47.3 (28.8)  & 71 (30) &  55.3 (31.3)\\
Illinois &  41.8 (26.4)  &  63.1 (25.7) &  50.7 (28.2) \\ \hline
\end{tabular}
\end{table}

\noindent Once again the difference in travel time is decreased if we consider the entire population. This is due to its structure. As there are more men and women who
use railroad in Illinois, the general average of travel time is increased.
This is reverse in New York.\\

\noindent In Table \ref{tab:label12} we test the equality of the averages when the ways of transportation to work is Bus/trolley bus.

\begin{table}
\caption{\label{tab:label12} Bus/trolleybus decision}
\centering
\begin{tabular}[h]{|c|c|c|} \hline 
 & Decision $n=1000$  & p-value\\ \hline \hline
Oracle test & not rejected & $0.19$\\
Expert test & not rejected & $0.13$ \\ \hline
Mixing test & not rejected & $0.75$ \\ 
\hline
\end{tabular}
\end{table}

\noindent In Table \ref{tab:label13} we reverse the set-up and we test the equality of the averages when the way of transportation to work is Walk.

\begin{table}
\caption{\label{tab:label13} Walk decision}
\centering
\begin{tabular}[h]{|c|c|c|} \hline 
 & Decision $n=1000$  &  p-value \\ \hline \hline
Oracle test & rejected & $0.05$ \\
Expert test & non-available & non-available \\ \hline
Mixing test & not rejected & $0.12$ \\ 
\hline
\end{tabular}
\end{table}

\section{Conclusion}\label{six}

\noindent From our point of view, one of the most interesting point is the usefulness of the varying mixing-weights model. It is a versatile model that can be used in many situations with missing microdata but aggregated information. The application treated exemplifies the modeling.\\
The second take-away message is the excellent performances of the Mixing test we propose. They can be guessed a priori thanks to the smallest eigenvalue of operators involved within the mixture-model. These nice performances were showed both theoretically and numerically.  \\
To conclude let us precise that this work can be easily extended to mixture-models with more than two components and can be done in a nonparametric setting when using the testing procedure proposed by Butucea and Tribouley \cite{BuTri06} as the Oracle test and the one given by Autin and Pouet \cite{AuPo} as the Mixing test.\\
An interesting extension that should really be considered is the case of mixing-weigths with errors.
This arises  when mixing-weigths are computed from a model with estimated parameters or from
experts' evaluation. In this case the solution of (\ref{algebra}) is no longer exact as the matrices $\Omega_X$ and $\Omega_Y$ are random. Preliminary simulation results tend to prove that moderate errors
have a small effect.

\section{Appendix}\label{seven}

\noindent In this section we provide the technical lemmas and the proposition required to prove the asymptotic normality under interest. For the sake of simplicity, we present them with respect to $X_1,\dots, X_n$ whereas an analogous version of them does exist for $Y_1,\dots, Y_n$. We recall that we assume that, for any $l \in \{1,2\}$ the mixing-weights of the model satisfy (\ref{ass-sup}).\\

\noindent Denote, for any $n \in \mathbb{N}^*$, any $l \in \{1,2\}$ and any $1 \leq i \leq n$ 
\begin{eqnarray}
W^{(l)}_{ni}&=&\frac{a_l(i)}{\sqrt{B^{(l)}_n}}\left(X_i-\omega_1(i)m_1-\omega_2(i)m_2\right) \label{W_i}.\label{W_ni}
\end{eqnarray}

\begin{lem}\label{lem-1}
For any $1 \leq i \leq n$
\begin{eqnarray*}
\mathbb{V}ar(X_i)&=&\sum_{l=1}^2 \left( \int \omega_l(i)(x-m_l)^2p_l(x)dx\right)+ \omega_1(i)\omega_2(i)(m_{1}+m_{2})^2.
\end{eqnarray*}
\end{lem}

\newpage
\noindent {\underline{Proof}}:\\

\noindent For any $1 \leq i \leq n,$
\begin{eqnarray*}
\mathbb{V}ar(X_i)&=&  \int (x-\omega_1(i)m_1+\omega_2(i)m_2)^2(\omega_1p_1(x)+\omega_2p_2(x))dx\\
&=&\sum_{l=1}^2 \left(\int \omega_l(i)(x-m_l)^2p_l(x)dx\right)+ \omega_1(i)\omega_2(i)(m_{1}+m_{2})^2.
\end{eqnarray*}

\noindent From Lemma \ref{lem-1}, we immediately derive:
\begin{lem}\label{lem-2}
For any $l \in \{1,2\}$, let $B_n^{(l)}$ be defined as in (\ref{B_n}). \begin{eqnarray*}
B_n^{(l)} & \geq& \min\left(\sigma_{1}^{2},\sigma_{2}^{2}\right)\sum_{i=1}^n a_{l}^2(i).
\end{eqnarray*}
\end{lem}

\noindent {\underline{Proof}}:\\
For any $l \in \{1,2\}$, by using Lemma \ref{lem-1} and the fact that, for any $1\leq i \leq n$ $\omega_1(i)+\omega_2(i)=1$,

\begin{eqnarray*}
B_n^{(l)} & \hspace{-0.3cm}=&\hspace{-0.2cm}  \sum_{i=1}^n a_{l}^{2}(i)  \mathbb Var(X_i)\\
& \hspace{-0.5cm}=& \hspace{-0.2cm}\sum_{i=1}^n \left[ a_{l}^{2}(i) \left(\sum_{l=1}^2 \int \omega_l(i)(x-m_l)^2p_l(x)dx\right)+ \omega_1(i)\omega_2(i)(m_{1}+m_{2})^2\right]\\
& \hspace{-0.5cm}\geq &\hspace{-0.2cm} \sum_{i=1}^n \hspace{-0.1cm} \left[a_{l}^{2}(i) \left(\sum_{l=1}^2  \omega_l(i)\right)\right] \hspace{-0.1cm} \min\left(\int (x-m_1)^2p_1(x)dx,\hspace{-0.1cm}\int (x-m_2)^2p_2(x)dx\hspace{-0.1cm}\right)\\
& \hspace{-0.5cm}=& \hspace{-0.2cm} \min\left(\sigma_{1}^{2},\sigma_{2}^{2}\right)\sum_{i=1}^n a_{l}^{2}(i).
\end{eqnarray*}

\begin{lem}\label{lem-3}
For any $l \in \{1,2\}$Let $B_n^{(l)}$ and $W^{(l)}_{ni}$ ($1 \leq i \leq n)$ be defined as in (\ref{B_n}) and (\ref{W_ni}). Then, for any $\epsilon >0$
\begin{eqnarray*}
\lim_{n\longrightarrow \infty} \ \sum_{i=1}^n \mathbb{E}\left[\left(W^{(l)}_{ni}\right)^2{\bf{1}}\left\{|W^{(l)}_{ni}|\geq \epsilon\right\}\right]= 0.
\end{eqnarray*}
\end{lem}

\noindent {\underline{Proof}}:\\

\noindent Fix $l \in \{1,2\}$. Let us define for any $n \in \mathbb{N}^*$
\begin{eqnarray*}
\kappa_n & = & \min\{ m_1,m_2\} + \frac {\varepsilon \sqrt{B^{(l)}_n}} {\sup_i\{|a_l(i)|\}}, \\
\kappa_n^\prime & = &  \max\{ m_1,m_2\} - \frac {\varepsilon \sqrt{B^{(l)}_n}} {\sup_i\{|a_l(i)|\}}.
\end{eqnarray*}

\begin{eqnarray*}
&&\sum_{i=1}^n \mathbb{E}\left[\left(W^{(l)}_{ni}\right)^2{\bf{1}}\left\{|W^{(l)}_{ni}|\geq \epsilon\right\}\right]\\
&=& \sum_{i=1}^n  \int_{|y| \geq \epsilon}
y^2  \ dF_{_{W^{(l)}_{ni}}}(x)\\
&\leq & \sum_{i=1}^n \frac {a^2_l(i)} {B^{(l)}_n} \int_{x>\kappa_n, x < \kappa_n^\prime} \hspace{-1cm}(x-\omega_1(i) m_1 - \omega_2(i) m_2)^2\ dF_{_{X_{i}}}(x)\\
&\leq & 2 \sum_{i=1}^n \frac {a^2_l(i)} {B^{(l)}_n} \left[ \sum_{l \in \{1,2\}}\omega_l(i)^3  \int_{x>\kappa_n, x < \kappa_n^\prime}
\hspace{-1cm}(x-m_l)^2 p_l(x) \ dx   \right. \\
& & \left.  \hspace{-0.2cm}\mathop{+}  \omega^2_2(i) \omega_1(i)  \int_{x>\kappa_n, x < \kappa_n^\prime}
\hspace{-1.3cm}(x-m_2)^2 p_1(x) \ dx +   \omega^2_1(i) \omega_2(i)  \int_{x>\kappa_n, x < \kappa_n^\prime}
\hspace{-1.3cm}(x-m_1)^2 p_2(x) \ dx \right]\\
&\leq & 2 \left(\sum_{i=1}^n \frac {a^2_l(i)} {B^{(l)}_n}\right) \ \sup_{n \in \mathbb{N}^*}\left(\sum_{(k,l) \in \{1,2\}^2}    \int_{x>\kappa_n, x < \kappa_n^\prime}
\hspace{-1cm}(x-m_k)^2 p_l(x) \ dx\right).\\
\end{eqnarray*}

\noindent Using Lemma \ref{lem-2}, for any $n \in \mathbb{N}^*$,
$$ \sum_{i=1}^n \frac {a^2_l(i)} {B^{(l)}_n}  \leq \left(\min(\sigma_1^2,\sigma_2^2)\right)^{-1}.$$

\noindent Then,  since variances under $p_1$ and $p_2$ are finite, the supremum over $n$  tends to $0$ in the integrals above, according to Lebesgue dominated convergence theorem.

\begin{lem}\label{lem-4} For any $1 \leq i \leq n$,
$$\mathbb{E}\left[(X_i-\mathbb{E}(X_i))^4\right]\leq C(m_1,m_2,p_1,p_2),$$
\noindent where $C(m_1,m_2,p_1,p_2):=32 \begin{displaystyle} \max_{(k,l) \in \{1,2\}^2} \int  \end{displaystyle} (x-m_k)^4p_l(x)dx.$ 
\end{lem}

\noindent {\underline{Proof}}:\\

\noindent We have
\begin{eqnarray*}
&&\mathbb{E}\left[(X_i-\mathbb{E}(X_i))^4\right] 
\\&=&\mathbb E\left[ \left( X_i - \omega_1(i) m_1 - \omega_2(i) m_2 \right)^4 \right] \\
&\leq & 8\
\omega_1(i)^4  \mathbb E\left( \left( X_i - m_1 \right)^4 \right) + 8\ \omega_2(i)^4 \mathbb E\left( \left( X_i - m_2 \right)^4 \right) \\ 
& = & 8  \left[ \omega_1(i)^5 \int (x-m_1)^4  p_1(x)\ dx    + \omega_1(i)^4 \omega_2(i) \int (x-m_1)^4  p_2(x)\ dx \right. \\
& &  \left. \: + \ \omega_2(i)^4 \omega_1(i) \int (x-m_2)^4 p_1(x)\  dx  + \omega_2(i)^5 \int (x-m_2)^4 p_2(x) \ dx \right]\\
& \leq &  32 \begin{displaystyle}\max_{(k,l) \in \{1,2\}^2} \int  \end{displaystyle} (x-m_k)^4p_l(x)dx.
\end{eqnarray*}

\begin{lem}\label{lem-5} For any $1 \leq l \leq 2$, the estimator 
$\begin{displaystyle}\hat{m}_l=\frac{1}{n}\sum_{i=1}^n \end{displaystyle}a_l(i)X_i$ of $m_l$ is consistent, that is $$\hat{m}_l\stackrel{Proba}{\to}  m_l.$$ 
\end{lem}

\noindent {\underline{Proof}}:\\

\noindent Let $\epsilon >0$ and $l\in \{1,2\}$. We have, using Bienayme-Chebyshev inequality and Lemma \ref{lem-1}

\begin{eqnarray*}
&&\mathbb P\left( \left| \hat m_l - m_l \right| > \varepsilon \right) \\
& \leq & \frac {1}{n^2 \epsilon^2} \sum_{i=1}^n a^2_l(i)
\mathbb{V}ar(X_i) \\
&= & \frac{1}{n^2 \varepsilon^2} \sum_{i=1}^n a^2_l(i)
\left[\sum_{l \in \{1,2\}} \left(\omega_l(i) \sigma^2_l +\omega_l(i)(1-\omega_l(i)) \ \frac{(m_1+m_2)^2}{2}\right)\right] \\
&\leq &  \frac {2}{n \varepsilon^2 K} \left[\sum_{l \in \{1,2\}} \sigma^2_l + \frac{(m_1+m_2)^2}{4}\right].\\
\end{eqnarray*}
Last inequality is obtained by using assumption on the
smallest eigenvalue of $\Omega'\Omega$, that is larger than $Kn$ (with $K>0$) and the fact that the supremum value for $x \in [0,1]$ of $x \to x(1-x)$ is equal to $\frac{1}{4}$. The right-hand side clearly tends to $0$ when $n$ goes to infinity.  We conclude that $\hat m_l$ is consistent.\\

\begin{prop}\label{prop-N} For any $l \in \{1,2\}$, any $n \in \mathbb{N}^*$ and any $1 \leq i \leq n$ consider $W^{(l)}_{ni}$, defined as in $(\ref{W_ni})$.
$$\sum_{i=1}^nW^{(l)}_{ni} \stackrel{\mathcal{L}}{\to}  \mathcal{N}(0,1).$$
\end{prop}

\noindent {\underline{Proof}}:\\

\noindent We apply Theorem 4.2 in Petrov \cite{Petrov}. It is the general setup for the Central Limit Theorem, for the  triangular array of series $(W_{ni})_{i,n}$ of independent random variables $X_i$ (that are not identically distributed).\\

\noindent 
If the three conditions are satisfied for any $\epsilon>0$ and any $\tau>0$
\begin{enumerate}
\item
$\begin{displaystyle}\lim_{n \to \infty}\sum_{i=1}^n \end{displaystyle} \mathbb P\left( |W^{(l)}_{ni} | \geq \varepsilon \right) = 0,$\\
\item 
$\begin{displaystyle}\lim_{n \to \infty} \sum_{i=1}^n \int_{|y| < \tau}\end{displaystyle}\ y dF_{_{W^{(l)}_{ni}}}(y)  =  0$,
\item $\begin{displaystyle}\lim_{n \to \infty}\sum_{i=1}^n \left\{ \int_{|y| < \tau} y^2 dF_{_{W^{(l)}_{ni}}}(y) - \left( \int_{|y| < \tau} y  dF_{_{W^{(l)}_{ni}}}(y) \right)^2 \right\} \end{displaystyle} =
1,$
\end{enumerate}
\noindent then \:
$ \begin{displaystyle}\sum_{i=1}^n W^{(l)}_{ni} \end{displaystyle}\stackrel{\mathcal{L}}{\to} \mathcal{N}(0,1).$\\

\noindent Let us prove that the three conditions are satisfied. Let $\epsilon>0$. Using Bienayme-Chebyshev inequality, 
\begin{eqnarray*}
\sum_{i=1}^n \mathbb P\left( | W^{(l)}_{ni} |  \geq   \varepsilon \right) 
& \leq &\epsilon^{-2}\sum_{i=1}^n \mathbb{E}\left[\left(W^{(l)}_{ni}\right)^2{\bf{1}}\left\{|W^{(l)}_{ni}|\geq \epsilon\right\}\right].
\end{eqnarray*}
\noindent Hence, condition $1$ is clearly satisfied by using Lemma \ref{lem-3}. \\

\noindent Let us move to condition $2$. We use the same trick as above. For any $\tau>0$
$$ \sum_{i=1}^n \int_{|y| < \tau} y dF_{_{W^{(l)}_{ni}}}(y)  =  \sum_{i=1}^n \left[ \int y dF_{_{W^{(l)}_{ni}}}(y) - \int_{|y|\geq \tau}
y dF_{_{W^{(l)}_{ni}}}(y)\right].$$
The first summand is equal to $0$ as the variables $W^{(l)}_{ni}$ are centered. 
\begin{eqnarray*} 
\sum_{i=1}^n \left| \int_{|y|\geq \tau} y dF_{_{W^{(l)}_{ni}}}(y) \right| &\leq & \sum_{i=1}^n \int_{|y|\geq \tau} |y| dF_{_{W^{(l)}_{ni}}}(y)\\
&\leq &\tau^{-1}\sum_{i=1}^n \mathbb{E}\left[\left(W^{(l)}_{ni}\right)^2{\bf{1}}\left\{|W^{(l)}_{ni}|\geq \tau\right\}\right].
\end{eqnarray*}
\noindent Condition $2$ is clearly satisfied by using Lemma \ref{lem-3}.\\

\noindent We end the proof with condition $3$. There are two parts (because of two summands) in this condition. For the first part  we proceed exactly as
in condition $2$. Indeed we have
$$ \sum_{i=1}^n \int_{|y| < \tau} y^2 dF_{_{W^{(l)}_{ni}}}(y) = \sum_{i=1}^n \int y^2 dF_{_{W^{(l)}_{ni}}}(y) - \sum_{i=1}^n
\int_{|y| \geq \tau} y^2 dF_{_{W^{(l)}_{ni}}}(y).$$
The first summand is exactly equal to $1$ and the second one tends to $0$ as $n$ goes to infinity, according to Lemma \ref{lem-3}.  Therefore it remains to prove that the second part tends to $0$ when $n$ goes to infinity. Because the variables $W^{(l)}_{ni}$ are centered and according to Cauchy-Schwarz inequality:
 \begin{eqnarray*}
 \sum_{i=1}^n \left( \int_{|y| < \tau} y dF_{_{W^{(l)}_{ni}}}(y) \right)^2 &=& \sum_{i=1}^n \left( \int_{|y| \geq \tau} y dF_{_{W^{(l)}_{ni}}}(y) \right)^2\\
 & \leq &  \sum_{i=1}^n \int_{|y| \geq \tau} y^2 dF_{_{W^{(l)}_{ni}}}(y) \\
&=& \sum_{i=1}^n \mathbb{E}\left[\left(W^{(l)}_{ni}\right)^2{\bf{1}}\left\{|W^{(l)}_{ni}|\geq \tau\right\}\right].
  \end{eqnarray*}

\noindent Still using Lemma \ref{lem-3}, we conclude that the second part we are interested in tends to $0$ when $n$ goes to infinity, as expected.\\

\end{document}